\documentclass[reqno,11pt]{amsart}

\usepackage{color}

\newtheorem{theorem}{Theorem}[section]
\newtheorem{lemma}[theorem]{Lemma}
\newtheorem{corollary}[theorem]{Corollary}

\theoremstyle{definition}
\newtheorem{assumption}[theorem]{Assumption}

\theoremstyle{remark}
\newtheorem{remark}[theorem]{Remark}
\newtheorem{example}[theorem]{Example}

\makeatletter 
 \def\dashint{%
 \operatorname%
 {\,\,\text{\bf--}\kern-.98em\DOTSI\intop\ilimits@\!\!}}
\def\dashnorm{\,\,\text{\bf--}\kern-.5em\|}
\def\ninf{\qopname\relax\@empty{inf\phantom{p}\!\!\!}}

 \makeatother

\newcommand\bC{\mathbb{C}}

\newcommand\bR{\mathbb{R}}
\newcommand\bS{\mathbb{S}}

\newcommand \bW{\mathbb{W}}

\newcommand \cL{\mathcal{L}}

\newcommand{\loc}{{\rm loc}\,}

\newcommand{\EO}{\overset{\scriptscriptstyle0}%
{E}\,\!}

 \newcommand{\mysection}[1]{\section{#1}
 \setcounter{equation}{0}}

\begin{document}

\title[Parabolic equations in Morrey spaces]
{On parabolic equations in Morrey spaces with
  VMO $a$ and Morrey $b,c$}
\author{N.V. Krylov}

\email{nkrylov@umn.edu}
\address{School of Mathematics, University of Minnesota, Minneapolis, MN, 55455}
 
 

\begin{abstract} 
We prove  existence and uniqueness
of solutions in Morrey spaces of functions 
with mixed norms
for
second-oder  parabolic equations  with
  VMO $a$ and Morrey $b,c$.
\end{abstract}

\maketitle

\mysection{Introduction}
                                                  \label{section 3.11.1}

Let $\bR^{d}$ be a $d-$dimensional Euclidean space of points
$x=(x^{1},...,x^{d})$ with $d\geq2$. Fix a $\delta\in(0,1]$ and let $\bS_{\delta}$ be the set
of $d\times d$ symmetric matrices whose eigenvalues
are in $[\delta,\delta^{-1}]$.
Assume that on $\bR^{d+1}=\{(t,x):t\in\bR,x\in\bR^{d}\}$ we are given a Borel real-valued $c$, an $\bR^{d}$-valued 
function  $b=(b^{i})$, and an $\bS_{\delta}$-valued $a
=(a^{ij})$. Define
$$
\cL u=\partial_{t}u+a^{ij}D_{ij}u+b^{i}D_{i}u,
$$
where the summation over the repeated indices
is enforced and we use the notation
$$
 D_{i}  =\frac{\partial}{\partial x^{i}} ,\quad Du= (D_{i}u),\quad
  D_{ij} = D_{i}D_{j}  ,\quad D^{2}u=(D_{ij}u),
\quad
 \partial_{t} =\frac{\partial}{\partial t} .
$$

We are going
to   investigate
the equation
\begin{equation}
                                         \label{6.15.2}
 \cL u-(c+\lambda)u=f
\end{equation}
 in $\bR^{d+1}$ and $[S,T]\times \bR^{d}$ with sufficiently large
$\lambda$. We are interested in rather regular
solutions and $b,c,f$ with rather poor integrability properties and, therefore, we work in Morrey spaces
$E^{1,2}_{p,q,\beta }$ (introduced in Section
\ref{section 3.27.1}).

G. Lieberman in \cite{Li_03} considered equations
in $E^{1,2}_{p,p,\beta }$ (no mixed norms)
in bounded domains. His assumptions on $b,c$
are much stronger than ours. For instance,
$|b|=1/|x|$ or $|c|=1/|x|^{2}$ (if $d\geq3$) do
not satisfy his integrability conditions.
In our setting to include those it suffices
to take $p=q>1$ close to one and, say $\beta_{0}=3$.
These functions are independent of time
and they do not fit in \cite{FHS_17}, \cite{CLV_96},
\cite{MT_13}
either, which
are about the elliptic equations with Morrey coefficients and solutions in Morrey class
in \cite{FHS_17} or in Sobolev class in \cite{CLV_96}
an \cite{MT_13}. In the case of elliptic equation
 one can find plenty of information in the recent article 
\cite{GO_22}. 

Much work has been done on elliptic equations in generalized Morrey spaces see, for instance,
\cite{PS_21}, where however there are no lower order terms. Many authors are investigating some rather exotic
problems like equations with measure data 
in \cite{Zh_19}, where again there are no lower order terms.
General information
about Morrey spaces and their generalizations can be found in \cite{No_19}, \cite{SHF_20}.

The literature on the parabolic equations in Morrey spaces
is not as rich as in the case of elliptic equations,
although it is worth drawing the reader's attention to
\cite{Li_03}, \cite{So_11}, \cite{Zh_19}, \cite{ZJZ_09}
and the reference therein.  Quite an interesting detail: In many papers the divergence form
equations are treated but in most of those with nondivergence form  
there are no lower order terms, and at the same time some of them
deal with the equations in   domains and claim that
having interior and boundary estimates 
in Morrey spaces provides global estimates
by usual method of partition of unity and flattening the boundary.
But this flattening necessarily introduces lower order terms,
treating which in Morrey spaces was first explicitly attempted rather late in \cite{FHS_17} and is  not an absolutely trivial task (see, for instance, the proof of the interpolation Lemma 3.9 in \cite{Kr_22} and compare it with the proof
of Theorem 4.3 of \cite{FHS_17} based on a false Lemma 4.2
of \cite{FHS_17}). It looks like this, looking absolutely
familiar, interpolation lemma for Morrey spaces
belongs to the author.

The author can indicate only two papers in which the lower order terms
in {\em parabolic\/} equations are in Morrey class 
(much wider than in \cite{Li_03}) and solutions
have two derivatives, these are \cite{Kr_22-2} and \cite{Kr_22-1}. 
As Hongjie Dong kindly pointed out (predicted) to the author,
Theorem 1.1 of \cite{Kr_22-1} follows from our
Theorem \ref{theorem 10.27.1} (which is new as far as the author is aware).

For the author the main motivation in developing
the theory in Morrey spaces
was to get the necessary tools to show that the corresponding It\^o stochastic equations admit solutions
for a larger than previously known class of $b$.
In \cite{RZ_20} such program was realized when
$a$ is the unit matrix and in the basic case
\begin{equation}
                                  \label{4.3.1}
 \Big(\int_{\bR}\Big(\int_{\bR^{d}}|b|^{p}\,dx\Big)^{q/p}
\,dt\Big)^{1/q}<\infty,\quad \frac{d}{p}+\frac{2}{q}=1.
\end{equation}
A very insightful observation in \cite{Ki_23} is that for the
solutions of It\^o stochastic equations to exist
in case of unit $a$ it suffices to have $s\in (1,d+2]$ and
$\rho_{0}>0$ such that for any standard parabolic cylinder $C_{r}$ with base being a ball of radius $r\leq \rho_{0}$ it holds that
\begin{equation}
                                  \label{4.3.02}
 \Big(r^{-d-2}\int_{C_{r}} |b|^{s}\,dx 
dt\Big)^{1/s}\leq \frak{b}r^{-1 },
\end{equation}
where $\frak{b}$ is sufficiently small.
Using H\"older's inequality one easily shows that
for $s=p\wedge q$ the left-hand side of \eqref{4.3.02}
is less than a constant times 
the left-hand side of \eqref{4.3.1} with 
$bI_{C_{r}}$ in place of $b$. The latter can be made
as small as we like on the account of choosing
$\rho_{0}$ sufficiently small.

The goal of this article is to combine, in a sense,
conditions \eqref{4.3.1} and \eqref{4.3.02} and develop
the solvability theory of \eqref{6.15.2}
in Morrey classes of functions without applying them
to It\^o stochastic equations.
This paper is a substantial expansion
of the PDE part of \cite{Kr_22-2}, where the
class of $b$ is more narrow and $c$ is bounded
but unconditional weak uniqueness of solutions
of It\^o's equation is achieved.
The expansion was prompted by the clever
observation in \cite{Ki_23} we were talking about
above.

Here is an example 
 in which our theory is applicable:
$|b|=ch$, where 
the constant $c>0$ is small enough and
$$
h(t,x)=I_{1>t>0,|x|<1}|x|^{-1}\Big(\frac{|x|}{\sqrt t}\Big)^{1/(d+1)},\,\,
a=2(\delta^{ij})+I_{x \ne0}\zeta(x )\sin(\ln|\ln |x |),
$$
where $\zeta$ is any smooth symmetric $d\times d$-matrix valued 
function vanishing for $|x|>1/2$ and
satisfying $|\zeta|\leq 1$.
This example is inadmissible in  
  \cite{RZ_20} because $b$ is too singular and
  $\sigma$ is not constant but admissible in
\cite{Ki_23} if $\sigma$ is constant
and yields weak solutions of It\^o's equations that are
conditionally weakly unique. The results of \cite{Kr_22-2} show that in this example solutions
are weakly unique without any additional conditions. 
By using the results of the present paper we show in
\cite{Kr_23} that any solution is, actually, strong and unique.

 The paper is organized as follows. 
In Section \ref{section 3.27.1} we introduce 
Morrey spaces and list some of their properties
with the proof of Lemma \ref{lemma 3.27.40}
deferred to Section \ref{section 4.2.1}.
Section \ref{section 3.27.3}
contains our main result
in parabolic and elliptic cases.  
Section \ref{section 4.2.2} is devoted to the study
of equations without lower order terms and
in Section \ref{section 4.2.3} we prove our
main Theorem \ref{theorem 10.27.1}

We conclude the introduction by some notation.
 
For $p,q\in[1,\infty)$ and domain $Q\subset\bR^{d+1}$ by $L_{p,q}(Q)$
we mean the space of Borel (real-, vector- or matrix-valued)
 functions on $Q$   with finite norm given {\em in
one of two ways which is fixed throughout the paper\/}:
\begin{equation}
                                 \label{3.27.3}
\|f\|_{L_{p,q}(Q)}^{q}=\|fI_{Q}\|_{L_{p,q}}^{q}
=\int_{\bR}\Big(\int_{\bR^{d}}|fI_{Q}(t,x)|^{p}\,
dx\Big)^{q/p}\,dt
\end{equation}
or
\begin{equation}
                                 \label{4.3.2}
\|f\|_{L_{p,q}(Q)}^{p}=\|fI_{Q}\|_{L_{p,q}}^{p}
=\int_{\bR^{d}}\Big(\int_{\bR}|fI_{Q}(t,x)|^{q}\,dt
 \Big)^{p/q}\,dx.
\end{equation}
It seems to the author that this paper is the first one where the norm \eqref{4.3.2} is systematically
used in the solvability theory of parabolic equations.

In case $f$ is a vector or matrix by $|f|$ we mean
the square root of the sum of squares of all entries
of $f$.
By the way, the norm from \eqref{4.3.2} is used
to treat $b$ from the above example.

Set $L_{p,q}=L_{p,q}(\bR^{d+1})$.
  As usual, we write something like $f\in L_{p,q,\loc}$
if $f\zeta\in L_{p,q}$ for any 
infinitely differentiable $\zeta$ with compact support. We write $\|u,v,..\|_{L_{p,q}}$ to mean
the sum of the $L_{p,q}$-norms of what is inside.

By $W^{1,2}_{p,q}(Q)$ we mean the collection
of $u$ such that $\partial_{t}u$, $D^{2}u $, $Du$, $u
\in L_{p,q}(Q)$. The norm in $W^{1,2}_{p,q}(Q)$
is introduced in an obvious way.
   We drop $Q$ if $Q=\bR^{d+1}$.

 Introduce
$$
B_{R}(x)=\{y\in\bR^{d}:|x-y|<R\},\quad B_{R} =B_{R}(0),
$$ 
$$
C_{\tau,\rho}(t,x)=[t,t+\tau)\times B_{\rho}(x),\quad C_{\rho}...=C_{\rho^{2},\rho}...,\quad C_{\rho}=C_{\rho}(0,0),
$$
and let $\bC_{\rho}$ be the collection of
$C_{\rho}(t,x)$.

If a Borel $\Gamma\subset \bR^{d+1}$, by $|\Gamma|$ we mean its Lebesgue
measure and
$$
\dashint_{\Gamma}f(x)\,dxdt:=\frac{1}{|\Gamma|}
\int_{\Gamma}f(x)\,dxdt.
$$
 If $C\in\bC_{\rho}$ we set
$$
\dashnorm f\|_{L_{p,q}(C)}=\|1\|_{L_{p,q}(C)}^{-1}
\|f\|_{L_{p,q}(C)}=N(d)\rho^{-d/p-2/q}\|f\|_{L_{p,q}(C)}.
$$

In the proofs of our results
we use various (finite) constants called $N$ which
may change from one occurrence to another
and depend on the data only in the same way
as it is
  indicated in the statements
of the results.

\mysection{Morrey spaces}
                             \label{section 3.27.1}

Take $p,q\in(1,\infty)$.
For    $\beta\geq 0$, introduce the
Morrey space $E_{p,q,\beta} $
as the set of $g\in  L_{p,q,\loc}$ such that  
\begin{equation}
                             \label{8.11.02}
\|g\|_{E_{p,q,\beta} }:=
\sup_{\rho\leq 1,C\in\bC_{\rho}}\rho^{\beta}
\dashnorm g  \|_{ L_{p,q}(C)} <\infty .
\end{equation}  
Define
$$
E^{1,2}_{p,q,\beta} =\{u:u,Du,D^{2}u,
\partial_{t}u\in E_{p,q,\beta} \}
$$
and provide $E^{1,2}_{p,q,\beta} $ with an obvious norm.
We are going to use some properties of these
spaces proved in \cite{Kr_22} when $E^{1,2}_{p,q,\beta}$
is introduced in terms of the $L_{p,q}$-spaces
with norm found in \eqref{3.27.3}.
The proofs extend to the case of the norm 
from \eqref{4.3.2} 
without any trouble. 
Also in \cite{Kr_22} the $E_{p,q,\beta}$-norm
is defined differently form \eqref{8.11.02}
without the restriction on $\rho$, just $\rho<\infty$.
In connection with this observe that if $\beta\leq
d/p+2/q$ and the support of $g$ is in some $C\in\bC_{1}$,
then the $E_{p,q,\beta}$-norm of $g$ is the same:
taken from \eqref{8.11.02} or from \cite{Kr_22}.
One more useful observation is that
\begin{equation}
                               \label{3.30.3}
\|g\|_{E_{p,q,\beta} }=\sup_{C\in \bC_{1}}
\|gI_{C}\|_{E_{p,q,\beta} }.
\end{equation}

\begin{remark}
                        \label{remark 3.27.1}
If $u$ vanishes outside $C\in\bC_{1}$, then, if $d/p +2/q \geq\beta $, one easily sees that
for any $r\geq 1$ and $C'\in \bC_{r}$
$$
r^{\beta }\dashnorm u\|_{L_{p,q}(C')}\leq N(d)
\|u\|_{L_{p,q}(C)}.
$$
\end{remark}

\begin{remark}
                       \label{remark 10.3.1}
A simple covering argument shows that
for $\beta >1$, $1\geq\rho\ge\rho'$ and $C\in\bC_{\rho}$
$$
\rho^{\beta}\dashnorm f\|_{L_{p ,q }(C)}\leq 
\nu(d,\beta,p,q,\rho')
\sup_{r\leq \rho'}r^{\beta}
\sup_{C\in \bC_{r}}\dashnorm f\|_{L_{p,q}(C)},
$$
where $\nu(d,\beta,p,q,\rho')=N(d)(\rho')^{d/p+2/q-\beta-d-2}$.

\end{remark} 

Here is a useful approximation result.
For functions $f(t,x)$ and $\varepsilon>0$ we define
$f^{(\varepsilon)}=f*\zeta_{\varepsilon}$, where
  $\zeta_{\varepsilon}=\varepsilon^{-d-1}\zeta
(t/\varepsilon,x/\varepsilon)$ with a nonnegative
$\zeta\in C^{\infty}_{0}(\bR^{d+1})$ which has unit integral and $\zeta(0)=1$.
\begin{lemma}
                               \label{lemma 3.14.3}
Let $f\in E_{p,q,\beta}$.
Then for any $C\in \bC$ and $\beta'>\beta $
\begin{equation} 
                            \label{3.14.10}
\lim_{\varepsilon\downarrow0}
\|(f-f^{(\varepsilon)})I_{C}\|_{E_{p,q,\beta'} }=0.
\end{equation}
\end{lemma}

Proof. Clearly, for any $r\in(0,1]$
$$
\|(f-f^{(\varepsilon)})I_{C}\|_{E_{p,q,\beta'} }
\leq 2r^{\beta'-\beta}\sup_{\rho\leq r}\rho^{\beta}
\sup_{C\in\bC_{\rho}}\dashnorm f\|_{L_{p,q}(C)}
$$
$$
+N(d,r,p,q) 
\|f-f^{(\varepsilon)}\|_{L_{p,q}(C)}.
$$
Here the first term on the right can be made as small as
we like uniformly in
$\varepsilon$ on the account of $r$ and the second
term tends to zero as $\varepsilon \downarrow 0$.
The lemma is proved. \qed

If $0<\beta\leq d/p+2/q$ we have the following as a consequence of the interpolation Lemma  5.10 of 
\cite{Kr_22}. 

\begin{lemma}
                        \label{lemma 3.27.4}
For any $u\in E^{1,2}_{p,q,\beta}$ and $\varepsilon\in(0,1]$,
\begin{equation}
                                    \label{3.20.6}
\| Du\|_{E_{p,q,\beta}}\leq \varepsilon
\|\partial_{t}u,D^{2}u\|_{E_{p,q,\beta}}
+N(d,p,q,\beta)\varepsilon^{-1}\| u\|_{E_{p,q,\beta}}.
\end{equation} 
\end{lemma}

If $\beta>d/p+2/q$, the result follows from
Theorem 10.2 of \cite{BIN_75},
since then
$$
\| f\|_{E_{p,q,\beta}}=N(d,p)\sup_{C\in\bC_{1}}
\|f\|_{L_{p,q}(C)}.
$$

If $u$ is independent of $t$, \eqref{3.20.6}
provides and interpolation inequality in the ``elliptic"
case, in which case for the first time in history it was stated as  
  Theorem 4.3 in \cite{FHS_17} but the proof is based on a false Lemma 4.2
of \cite{FHS_17}.
 
The following embedding lemma is proved in Section \ref{section 4.2.1}.
 
\begin{lemma}
                        \label{lemma 3.27.40}
(A) Let (a) $ \beta <2 $
or (b) $d/p+2/q<2$. Then any $u\in E^{1,2}_{p,q,\beta}$
is bounded and continuous  and for any $\varepsilon
\in(0,1]$ in case (a)
\begin{equation}
                                    \label{3.20.06}
|u|
\leq \varepsilon^{2-\beta} \|\partial_{t}u,D^{2}u\|_{E_{p,q,\beta}}+N(d, \beta)\varepsilon^{-\beta} \| u\|_{E _{p,q,\beta}}
\end{equation} 
with $\beta=2-d/p-2/q$ in the powers of $\varepsilon$ 
and $N(d,p,q)$ in \eqref{3.20.06} in case (b).

(B) Let (a) $  \beta<1 $
or (b) $d/p+2/q<1$. Then, for any $u\in E^{1,2}_{p,q,\beta}$, $Du$
is bounded and continuous  and for any $\varepsilon
\in(0,1]$ in case (a)
\begin{equation}
                                    \label{3.20.60}
|Du|
\leq \varepsilon^{1-\beta} \|\partial_{t}u,D^{2}u\|_{E_{p,q,\beta}}+N(d, \beta)\varepsilon^{-\beta-1} \| u\|_{E _{p,q,\beta}}
\end{equation} 
with $\beta=1-d/p-2/q$ in the powers of $\varepsilon$ and $N(d,p,q)$  in \eqref{3.20.60} in case (b).
\end{lemma}

Another useful property of $E^{1,2}_{p,q,\beta}$
is an embedding theorem, which
follows from Corollary 5.7 of \cite{Kr_22}.

\begin{lemma}
                        \label{lemma 3.27.5}
Let $ 1< \beta\leq d/p+2/q$. Then for any $u\in E^{1,2}_{p,q,\beta}$
we have
\begin{equation}
                                    \label{3.27.5}
\|Du\|_{E_{r,s,\beta-1}}\leq N(d,p,q,\beta)  \| u\|_{E^{1,2}_{p,q,\beta}},
\end{equation} 
where $(r,s)(\beta-1)=(p,q)\beta$.
\end{lemma}

One more embedding result is 
a consequence of Theorem 5.6 of \cite{Kr_22}
with $\alpha=2$.

\begin{lemma}
                           \label{lemma 4.1.1}
Let $2 <\beta\leq d/p+2/q$, $(r,s)(\beta-2)= 
(p,q)\beta$, $u\in E^{1,2}_{p,q,\beta}$. Then 
 there exists $N=N(d,p,q,\beta)$ such that 
 
\begin{equation}
                                 \label{4.1.2}
\|u\|_{E_{r,s,\beta-2}}\leq N(d,p,q,\beta)  \| u\|_{E^{1,2}_{p,q,\beta}}.
\end{equation}

\end{lemma}

The following  complementary result not used
in this paper is proved in Section \ref{section 4.2.1}.

\begin{lemma}
                             \label{lemmq 4.6.1}
Let $1<\beta\leq d/p+2/q$, $\beta<2$. Then for any
$u\in E^{2}_{p,q,\beta}$, $\rho\leq 1$,
$(t_{i},x_{i})\in C_{\rho}$, $i=1,2$, we have
\begin{equation}
                           \label{4.6.1}
|u(t_{1},x_{1})-u(t_{2},x_{2})|\leq
N(d,p,q,\beta)\rho^{2-\beta}\|u\|_{E^{1,2}_{p,q,\beta}}.
\end{equation}

In case $0<\beta\leq d/p+2/q$, $\beta<1$, we have
\begin{equation}
                           \label{4.6.2}
|Du(t_{1},x_{1})-Du(t_{2},x_{2})|\leq
N(d, \beta)\rho^{1-\beta}\|u\|_{E^{1,2}_{p,q,\beta}}.
\end{equation}

\end{lemma}

Observe that, if we want to have functions from
$W^{1,2}_{p,q}$ to be H\"older continuous
(as in \eqref{4.6.1}), we need $d/p+2/q<2$.
For H\"older continuity of the gradients of
$W^{1,2}_{p,q}$-functions we need $d/p+2/q<1$.
There are no restrictions on $d/p+2/q$ from above
in Lemma \ref{lemmq 4.6.1}.

\mysection{Main result}
                            \label{section 3.27.3}

We have  $p ,q ,\beta_{0}   $ such that
\begin{equation}
                        \label{3.21.01}
p ,q  \in(1,\infty),\quad\beta_{0}\in(0,\infty), \quad \beta_{0}\ne 1, 2,
\quad \frac{d}{p }+\frac{2}{q }\geq\beta_{0}.
\end{equation}

Fix some $\rho_{a},\rho_{b,c} \in(0,1]$.
Parameters $\theta$, $\check b$  and $\check c$  below will be specified later.
\begin{assumption}
                        \label{assumption 12.12.3}
 We have   
\begin{equation}
                                 \label{6.3.1}
 a^{\sharp}_{x,\rho_{a}}:=\sup_{\substack{\rho\leq\rho_{a}\\C\in\bC_{\rho}}}\dashint_{C}|a (t,x)-a_{C}(t)|
\,dx dt \leq  \theta,
\end{equation}
where
$$
a_{C}(t)=\dashint_{C}a (t,x)\,dxds\quad (\text{note}\quad t\quad \text{and}\quad ds).
$$
\end{assumption}

\begin{assumption}
                      \label{assumption 3.14.1}
If $\beta_{0}<1$,
\begin{equation}
                           \label{4.1.5}
\mathfrak{b} :=\rho_{b,c}^{1-\beta_{0}}\sup_{r\leq \rho_{b,c}}r^{\beta_{0}}
\sup_{C\in \bC_{r}} 
\dashnorm b \|_{L_{p ,q }(C)}<\infty .
\end{equation}
If $\beta_{0}>1$, then  
\begin{equation}
                           \label{3.14.2}
\mathfrak{b}_{\rho_{b,c}}:=\sup_{r\leq\rho_{b,c}}r
\sup_{C\in \bC_{r}} 
\dashnorm b \|_{L_{p\beta_{0},q\beta_{0}}(C)}\leq   \check b .
\end{equation}
\end{assumption}

\begin{assumption}
                      \label{assumption 3.27.2}
If $\beta_{0}<2$,
\begin{equation}
                           \label{4.1.6}
\mathfrak{c}:=\rho_{b,c}^{2-\beta_{0}}\sup_{r\leq \rho_{b,c}}r^{\beta_{0}}
\sup_{C\in \bC_{r}} 
\dashnorm c \|_{L_{p ,q }(C)}<\infty .
\end{equation}
If $\beta_{0}> 2$,
\begin{equation}
                           \label{4.2.1}
\mathfrak{c}_{\rho_{b,c}} :=\sup_{r\leq \rho_{b,c}}r^{2}
\sup_{C\in \bC_{r}} 
\dashnorm c \|_{L_{p\beta_{0},q\beta_{0}}(C)}
\leq \check c.
\end{equation}
 
\end{assumption}

Let us specify $\theta$ in \eqref{6.3.1}.   It is easy to choose
$\theta_{1}(d,\delta,p,q)$ introduced in Lemma \ref{lemma 3.25.1}, so that it is a decreasing
function of $d$, and we suppose it is done.
In the following $\alpha\in(0,1)$ is a free parameter.

\begin{assumption}
                    \label{assumption 3.25.1}
For   $r$ defined as the least number such that
$$
 r\geq  (d+2)/\alpha,\quad r\geq p,q 
$$
and $\Theta(\alpha)=\{(p',q'):
p\leq p'\leq r,q\leq q'\leq r\}$
Assumption \ref{assumption 12.12.3} is satisfied with 
$$
\theta=\inf_{\Theta(\alpha)} \theta_{1}(d+1,\delta,p',q')  =: 
\theta (d,\delta,p,q,\alpha) .
$$
\end{assumption}

The fact that this $\theta>0$  follows easily
from the arguments in \cite{DK_18} leading to defining $\theta_{1}(d,\delta,p,q)$ in Lemma \ref{lemma 3.25.1}.

Recall that
$$
\cL u(t,x)=\partial_{t}u(t,x)+ a^{ij}(t,x)D_{ij}u(t,x)+b^{i}(t,x)D_{i}u(t,x).
$$
 
Here is our main result. 
In the following theorem the constants $\check \lambda_{0}$
and $N$ depend on $\mathfrak{b}$ only if $\beta_{0}< 1$
and depend on $\mathfrak{c}$ only if $\beta_{0}<2$, but in order not to make the statement too heavy
we list these parameters anyway.
The constant $\nu=\nu(d,\beta_{0},p,q,\rho_{b,c})$
is taken from Remark \ref{remark 10.3.1}.
\begin{theorem}
                       \label{theorem 10.27.1}
 
Under the above assumptions there exist 
$$
(\check b,\check c) =(\check b,\check c)(d,\delta, p, 
q, \rho_{a}, \beta_{0},\alpha )\in(0,1]^{2},
$$
$$
\check \lambda_{0}=\check \lambda_{0}
(d,\delta, p, 
q, \rho_{a}, \beta_{0},\alpha,\mathfrak{b} ,\mathfrak{c})>0
$$
 such that,
if \eqref{3.14.2} holds with this $ \check b$
in case $\beta_{0}>1$ and \eqref{4.2.1} holds with this $\check c$ in case $\beta_{0}>2$,  
then for  any
$\lambda\geq  \rho_{b,c}^{-1}\check \lambda_{0} $  and $f\in E_{p,q,\beta_{0} }$ there exists
a unique $E^{1,2}_{p,q,\beta_{0} }$-solution $u$
of $\cL u-(c+\lambda) u=f$. Furthermore, there exists
a constant $N$ depending only on  $d$, $\delta$, $p$, $q$, $\rho_{a}$,  $\beta_{0}$,
$\alpha$, $\mathfrak{b}$, $\mathfrak{c}$,    
such that 
\begin{equation}
                                \label{10.28.01}
\|\partial_{t} u,D^{2}u, \sqrt\lambda Du, \lambda u\|_{E_{p,q,\beta_{0} }}
\leq   N \nu^{-1}\|f\|_{E_{p,q,\beta_{0} }}.    
\end{equation}
\end{theorem}

\begin{remark}
                     \label{remark 8.1.1}
If $-\infty\leq S<T\leq \infty$,
Theorem \ref{theorem 10.27.1} also holds if we replace
$E_{p,q,\beta_{0} }$ with $E_{p,q,\beta_{0} }((S,T)\times\bR^{d})$ defined as the set of functions $f$ on $(S,T)\times\bR^{d}$
such that $fI_{(S,T)\times\bR^{d}}\in
E_{p,q,\beta_{0} }$. 
Of course, we also
replace $E^{1,2}_{p,q,\beta_{0}}$ with
$\EO^{1,2}_{p,q,\beta_{0} }((S,T)\times\bR^{d})$,
which is the set of $u$ such that $(u,Du,D^{2}u,
\partial_{t}u)I_{(S,\infty)\times\bR^{d}}\in
E_{p,q,\beta_{0} }$ and $u(t,\cdot)=0$ for $t>T$.

Indeed, in case $T=\infty$, the existence of solutions and the estimate is trivial.
In what concerns uniqueness, having in mind
the possibility to extend functions given 
on $(S,\infty)$ 
to $\bR$, it suffices to prove that, 
if 
$v\in E^{1,2}_{p,q,\beta_{0}}$ is such that
$\cL v-(c+\lambda) v=0$ for $t>S$, then $v=0$ for $t>S$. This is an easy task because solutions
of $\cL v-(c+\lambda) v=g$ can be obtained by the method of continuity and, if $\cL=\partial_{t}+\Delta$ and $c=0$,
this property is well known and follows from
the explicit formula for the solutions.
The case of general $T<\infty$, follows easily from what is said above (just set $f(t,\cdot)=0$ for
$t>T$).

\end{remark}

\begin{remark}
                     \label{remark 10.27.2}   
The unique solution $u$ from Theorem \ref{theorem 10.27.1}
possesses the following properties

a) obviously, $u\in W^{1,2}_{p,q,\loc}$;

b) for $\beta_{0}>1$, by Lemma \ref{lemma 3.27.5}, we have
$Du\in L_{r,s,\loc} $, where $(r,s)=  (\beta_{0}-1)^{-1}\beta_{0}(p ,q )$;

c) if $\beta_{0}<2$, $u$ is bounded and (H\"older) continuous and, if $\beta_{0}<1$, then $Du$ is bounded and (H\"older)  continuous.

\end{remark}

\begin{example}
                          \label{example 5.12.1}
It turns out that $\check b$ is generally small.
To show this take $d>4$ and set
$$
p=q=\frac{d-1}{2},\quad \beta_{0}= \frac{2d-1}{d-1},
\quad p'=\frac{d+2}{\beta_{0}},
$$
so that  
$ 2<\beta_{0}< d/p+2/q $, 
$\beta_{0}p<d$, $p'>p$. Then take any
 $\lambda>0$, and consider the
function
$$
u(t,x)=e^{\lambda t+|x|^{2}/(4t)}I_{t<0}.
$$
This function satisfies
$$
\partial_{t}u+\Delta u-\frac{d}{|x|^{2}}x^{i}D_{i}u
-\lambda u=0
$$
in $\bR^{d+1}$. Also, as is easy to check, $p'<(d+2)/2$ and
$v:=u,Du, D^{2}u,\partial_{t}u\in L_{p',p'}$ and, hence,
for any $C\in \bC_{\rho}$
$$
\rho^{\beta_{0}}\dashnorm v\|_{L_{p,p}(C)}
\leq \rho^{\beta_{0}}\dashnorm v\|_{L_{p',p'}(C)}
=N(d)\|v\|_{L_{p',p'}(C)}\leq N(d)\|v\|_{L_{p',p'}},
$$
so that $u\in E^{1,2}_{p,q,\beta_{0}}$.
 In addition, for $b=-xd/|x|^{2}$ the quantity
$\mathfrak{b}_{\rho_{b,c}}$ is finite (and independent of
$\rho_{b,c}$) since $p\beta_{0}<d$.
However, estimate \eqref{10.28.01} fails,
which shows that the constant factor $d$ in $b$
is not sufficiently small.

\end{example}

Theorem \ref{theorem 10.27.1} is applicable
in the case of $a,b,c,f$   independent of $t$,
and owing to uniqueness yields solutions that
are also independent of $t$. This leads to
the ``elliptic'' counterpart, in which
$E_{p,\beta_{0}}$ and $E^{1,2}_{p,\beta_{0}}$
are subsets of $E_{p,q,\beta_{0}}$ and $E^{1,2}_{p,q,\beta_{0}}$ consisting of time-independent functions
and which says that in that case Theorem
\ref{theorem 10.27.1} holds true if we replace
$E_{p,q,\beta_{0}}$ and $E^{1,2}_{p,q,\beta_{0}}$ 
with $E_{p,\beta_{0}}$ and $E^{1,2}_{p,\beta_{0}}$
and $\cL u$ with $a^{ij}D_{ij}u+b^{i}D_{i}u+cu$. 
 
\begin{remark}
                            \label{remark 4.2.1}
Make the change $u(t,x)=v(t/\rho_{b,c}^{2},
x/\rho_{b,c})$ of the unknown function $u$ in
  the equation $\cL u-(c+\lambda)u=f$. Then $v$
will satisfy
$$
\partial_{t}v+\tilde a^{ij}D_{ij}v+\tilde b^{i}D_{i}v
-(\tilde c+\tilde \lambda)v=\tilde f,
$$
where $\tilde \lambda=\rho_{b,c}^{2}\lambda$,
$$
(\tilde a,\tilde b,\tilde c,\tilde f)(t,x)
=(a,\rho_{b,c}b,\rho_{b,c}^{2}c,\rho_{b,c}^{2}f)
(\rho_{b,c}^{2}t,\rho_{b,c}x).
$$
It is easy to check that, since $\rho_{b,c}\leq 1$, $\tilde a,\tilde b,\tilde c$
satisfy Assumptions \ref{assumption 12.12.3}
and \ref{assumption 3.14.1} with {\em the same\/} $\rho_{a}$,
$\theta$, $\mathfrak{b}$, $\check b$, $\mathfrak{c}$, 
$\check c$ {\em and\/} 1 in place of $\rho_{b,c}$.
Then assuming that Theorem \ref{theorem 10.27.1}
is true for $\rho_{b,c}=1$ we conclude that
for $\tilde \lambda\geq \check\lambda_{0}$
\begin{equation}
                                \label{4.2.3}
\|\partial_{t} v,D^{2}v, \sqrt{\tilde\lambda} Dv, \tilde\lambda v\|_{E_{p,q,\beta_{0} }}
\leq   N  \|\tilde f\|_{E_{p,q,\beta_{0} }}.    
\end{equation}
Since $\rho_{b,c}\leq 1$, $\|\tilde f\|_{E_{p,q,\beta_{0} }} \leq \rho_{b,c}^{2}\| f\|_{E_{p,q,\beta_{0} }} $. Also
$$
\|\partial_{t} v,D^{2}v, \sqrt\lambda Dv, \lambda u\|_{E_{p,q,\beta_{0} }}=\rho_{b,c}^{2}\|\partial_{t} u,D^{2}u, \sqrt{\tilde\lambda/\rho^{2}_{b,c}} Du,  u\tilde\lambda/\rho^{2}_{b,c}\|_{E'_{p,q,\beta_{0} }},
$$
where the norm in $E'_{p,q,\beta_{0} }$ is defined
as in $E_{p,q,\beta_{0} }$ but with $\rho\leq 1$
in \eqref{8.11.02}
being replaced with $\rho\leq \rho_{b,c}$.
Now 
Remark \ref{remark 10.3.1} allows us to obtain
\eqref{10.28.01} from \eqref{4.2.3}.
\end{remark}

This remark shows that without loosing generality
we impose

\begin{assumption}
                           \label{assumption 4.3.1}
We have $\rho_{b,c}=1$.
\end{assumption}

\mysection{Model case}

                           \label{section 4.2.2}

Here we take $p,q\in(1,\infty)$ not necessarily related to \eqref{3.21.01}.
Introduce
$$
\cL_{0} u(t,x)=\partial_{t}u+ a^{ij}(t,x)D_{ij}u(t,x).
$$

Here is a corollary of Theorem 6.3
of \cite{DK_18}.

\begin{theorem}
                       \label{theorem 12.18.1}
There exists $\theta =\theta_{0} (d,\delta,p,q)>0$
such that with that $\theta$  in Assumption \ref{assumption 12.12.3},   there exists   constants
$$
(\lambda_{0},N_{0})=(\lambda_{0},N_{0})(d,\delta,p,q, \rho_{a}),\quad \lambda_{0}>0,
$$
such that for any $\lambda\geq \lambda_{0}$  and any
$f\in L_{p,q} $  there exists a unique
$u\in W^{1,2}_{p,q} $ such that
  $\cL_{0} u -\lambda u =f $. 
  Moreover, 
\begin{equation}
                                 \label{3.27.8}
\|\lambda u,\sqrt\lambda Du, D^{2}u\|_{L_{p,q} }\leq N_{0}
\|f\|_{L_{p,q} }.
\end{equation}

\end{theorem}

As in Remark 3.3 of \cite{Kr_pt1} we note that
Theorem \ref{theorem 12.18.1} is proved in 
\cite{DK_18} only  as an a priori estimate
in spaces generated by the norm from \eqref{3.27.3}
but yields the solvability in a usual way
by the method of continuity
and that we can replace \eqref{3.27.3}
with \eqref{4.3.2} (that has been already mentioned in 
\cite{DK_18}).  
 
If $Q\subset \bR^{d+1}$ by $\partial' Q$
we mean the set of points $(t ,x )\in \partial
Q$, for each of which there is $\tau>0$ and
a continuous curve $(s,x_{s})$, $s\in [t -\tau,t ]$,
such that $x_{t}=x$ and $(s,x_{s})\in Q$, $s\in [t -\tau,t )$.
  
\begin{corollary}
                            \label{corollary 3.18.1}

Assume that $a$ is smooth. Take
$R\leq 1$, smooth $f$, and let $u$ be the classical solution of
$$
\cL_{0} u+f=0
$$
in $C_{R}$ with zero boundary condition
on $\partial' C_{R}$. Then
$$
\dashnorm u\|_{L_{p,q}(C_{R})}\leq NR^{2 }\dashnorm  f\|_{L_{p,q}(C_{R})},
$$
where $N$ depends only on $d,\delta,p,q, \rho_{a}$.
\end{corollary}

Indeed, the case of general $R\leq 1$ is reduced to
$R=1$ by using parabolic dilations. If $R=1$, 
the maximum principle allows us to concentrate on $f\geq0$ and also shows that
$u(t,x)e^{\lambda_{0}t}$ is smaller  in $C_{1}$
than the solution $v$ of
$$
\cL v-\lambda_{0} v+I_{C_{1}}fe^{\lambda_{0}t}=0
$$
in $\bR^{d+1}$. Then our assertion follows from \eqref{3.27.8}.

\begin{corollary}
                     \label{corollary 10.4.01}
There exists $N $,
depending only on $d,\delta,p,q, \rho_{a}$, such that
for any  
$u\in  W^{1,2}_{p ,q }$ we have
\begin{equation}
                            \label{10.2.23}
\|D^{2}u,Du,u\|_{L_{p,q}}
\leq N \|\cL_{0} u\|_{L_{p,q}}+N\|  u\|_{L_{p,q}}.
\end{equation}
\end{corollary}

By using a standard way of localizing
  such estimates given, for instance,
in the proof of Lemma 2.4.4 of \cite{Kr_08},
one arrives at the following.
\begin{lemma}
                        \label{lemma 10.4.3}
Suppose that Assumption \ref{assumption 12.12.3}
is satisfied with $\theta_{0}(d,\delta,p,q)$ from Theorem \ref{theorem 12.18.1}.
Then there exists $N $,
depending only on $d,\delta,p,q, \rho_{a}$, such that for any $0<\rho_{1}<\rho_{2}\leq 1$ and
$u\in  W^{1,2}_{p ,q,\loc }(C_{\rho_{2}})
\cap C(\bar C_{\rho_{2}})$ we have
\begin{equation}
                            \label{10.4.4}
\|D^{2}u\|_{L_{p,q}(C_{\rho_{1}})} 
\leq N \|\cL_{0} u\|_{L_{p,q}(C_{\rho_{2}})}+N
(\rho_{2}-\rho_{1})^{-2}\|  u-l\|_{L_{p,q}(C_{\rho_{2}})},
\end{equation}
\begin{equation}
                            \label{3.25.3}
\|Du\|_{L_{p,q}(C_{\rho_{1}})} 
\leq N \|\cL_{0} u\|_{L_{p,q}(C_{\rho_{2}})}+N[\rho_{1}^{-1}+
(\rho_{2}-\rho_{1})^{-2}]\|  u \|_{L_{p,q}(C_{\rho_{2}})},
\end{equation}
where $l$ is any  affine function of $x$.
\end{lemma}

Here the appearance of $l$ is explained by the fact that $D^{2}u$ and $\cL_{0} u$ do not change
if we take $u-l$ in place of $u$ and \eqref{3.25.3}
follows from \eqref{10.4.4} and  interpolation inequalities.

One finds in \cite{DK_18} a discussion of
initial-boundary problems, which leads to the following.
\begin{lemma}
                         \label{lemma 3.25.1}
Let  $f\in L_{p,q}(C_{1})$.
There exists $\theta_{1}=\theta_{1}(d,\delta,p,q)>0$
($\theta_{1}\leq \theta_{0}$) such that under  Assumption \ref{assumption 12.12.3} with
 $\theta=\theta_{1}(d,\delta,p,q)$
 there is a unique $u\in W^{1,2}_{p,q}(C_{1})$
vanishing on $\partial' C_{1}$ such that
$\cL_{0}u=f$ in $C_{1}$. Furthermore,
$$
\|u\|_{W^{1,2}_{p,q}(C_{1})}\leq N_{1}\|f\|_{L_{p,q}(C_{1})},
$$
where $N_{1}=N_{1}(d,\delta,p,q,\rho_{a})$.
\end{lemma}

Starting from this point we suppose that  
Assumption \ref{assumption 3.25.1} is satisfied
and thus the restriction on $\theta$ is finally fixed.

\begin{corollary}
                       \label{corollary 3.25.1}
Let 
$v\in W^{1,2}_{p,q }(C_{\rho})$, for any $\rho<1$, $\cL_{0}v=0$
in $C_{1}$.  
Then, for any $\rho<1$, we have $v\in W^{1,2}_{r,r }(C_{\rho})$ and
\begin{equation}
                           \label{3.25.1}
\|v\|_{W^{1,2}_{r,r }(C_{\rho})}\leq
N \|v\|_{L_{p,q }(C_{1})},
\end{equation}
\begin{equation}
                           \label{3.26.1}
\sup_{C_{\rho}}|v| \leq
N \|v\|_{L_{p,q }(C_{1})},
\end{equation}
where the constants $N$ depend  only on $\alpha,\rho, d,\delta,\rho_{a},p,q$.
\end{corollary}

Indeed, estimate \eqref{3.26.1}
follows from \eqref{3.25.1} by embedding theorems
($r>d+2$). To prove
\eqref{3.25.1} we use a quite standard
 argument  
(see, for instance, the proof of Theorem 4.3.12
of \cite{Kr_08}).
Take a $\rho'<1$ and a $\zeta\in C^{\infty}_{0}(\bR^{d+1})$
vanishing in $C_{(1+\rho')/2}\setminus C_{\rho'}$
and equal to one on $C_{\rho'}$. Then
$\cL_{0}(\zeta v)$ only contains $v$ and $Dv$,
which are in $L_{p',q'}(C_{(1+\rho')/2})$
by embedding theorems (see 
Theorem 10.2 in \cite{BIN_75}),
provided
$$
\frac{2}{q}+\frac{d}{p}\leq 1+\frac{2}{q' }+\frac{d}{p '}.
$$
We may take $p\leq p'\leq r,q\leq q'\leq r$ and then
by Lemma \ref{lemma 3.25.1} we have $v\in
W^{1,2}_{p',q'}(C_{ \rho' })$. The arbitrariness
of $\rho'$ shows that 
$v\in
W^{1,2}_{p',q'}(C_{ \rho})$ for any $\rho<1$. Then we can do this procedure again and again and in finitely
many
 steps we will come to the conclusion that
$v\in W^{1,2}_{r,r }(C_{\rho})$ for any $\rho<1$.
By accompanying these steps with the corresponding estimates we see that
$$
\|v\|_{W^{1,2}_{r,r }(C_{\rho})}\leq
N \|v\|_{W^{1,2}_{p,q }(C_{(1+\rho)/2})},
$$
and it only remains to use Lemma \ref{lemma 10.4.3}
and interpolation inequalities
to get \eqref{3.25.1}.

\begin{corollary}
                       \label{corollary 3.25.2}
Let  
$\rho\leq 1$,
$v\in W^{1,2}_{r,r }(C_{\rho})$, and $\cL_{0}v=0$
in $C_{\rho}$. Then   for any $\rho<1$ we have
\begin{equation}
                           \label{3.25.01}
\dashnorm D^{2}v\|_{L_{r,r}(C_{ \rho/2})}\leq
N  \rho^{-2}\dashnorm v-l\|_{L_{p,q }(C_{\rho})}
\end{equation}
where $N$ depends only on $  d,\delta,\rho_{a},p,q
,\alpha$ (independent of $\rho$) and $l$ is an arbitrary
affine function.
\end{corollary}

Indeed, parabolic scaling reduce the whole
situation to the one in which $\rho=1$,
and then the result follows from \eqref{3.25.1}.

\begin{lemma}
                         \label{lemma 10.4.1}
Let $p,q\in(1,\infty)$. Then there is  $N,N_{1},N_{2}$ depending
only on $d,\delta,p,q,\alpha$, and $ \rho_{a}$ such that,
if  
$0<4\rho_{1}\leq \rho_{2}\leq 1$,  and
$u\in  W^{1,2}_{p ,q  }(C_{\rho_{2}})$, then
\begin{equation} 
                              \label{7.15.3}
\dashnorm D^{2}u\|_{L_{p,q}(C_{\rho_{1}})}
\leq N_{1} (\rho_{2}/\rho_{1})^{\gamma}\dashnorm \cL_{0} u\|_{L_{p,q}(C_{\rho_{2}})}
+N  (\rho_{2}/\rho_{1})^{\alpha}\rho_{2}^{-2}\dashnorm u-l\|_{L_{p,q}(C_{\rho_{2}})},
\end{equation} 
\begin{equation}
                                  \label{7.15.30}
\dashnorm D^{2}u\|_{L_{p,q}(C_{\rho_{1}})}
\leq N_{1} (\rho_{2}/\rho_{1})^{\gamma}\dashnorm \cL_{0} u\|_{L_{p,q}(C_{\rho_{2}})}
+N_{2} (\rho_{2}/\rho_{1})^{\alpha}\dashnorm D^{2}u\|_{L_{p,q}(C_{\rho_{2}})}   , 
\end{equation}
where  $l$ is any  affine function of $x$ and
$\gamma=2+d/p+2/q$.
\end{lemma}

Proof.   Having in mind obvious approximations we may assume that $u$ and $a$ are  infinitely differentiable
in $\bar C_{\rho_{2}}$. Then by classical results there is smooth solution $v$ of the
equation $\cL_{0}v=0$ in $C_{\rho_{2}}$ with boundary
data $v=u$ on $\partial'C_{\rho_{2}}$. This
solution belongs to $\bW^{1,2}_{r,r}(C_{\rho_{2}})$ for any $r\geq 1$. We take $r$ from Assumption 
\ref{assumption 3.25.1}.   
 Then by Corollary \ref{corollary 3.25.2}  we have
$$  
\dashnorm D^{2}v\|_{L_{r,r}(C_{ \rho_{2}/2})}
\leq N  \rho_{2} ^{-2}\dashnorm v-l\|_{L_{p,q}(C_{ \rho_{2}} )}
$$
$$
\leq N  \rho_{2} ^{-2}\dashnorm u-l\|_{L_{p,q}(C_{ \rho_{2}} )}+N  \rho_{2} ^{-2}\dashnorm w\|_{L_{p,q}(C_{ \rho_{2}} )},
$$
where $l$ is any affine function and 
$w=u-v$. Since $w=0$ on $\partial'    
C_{\rho_{2}}$ and $\cL_{0}w=\cL_{0} u $, by Corollary \ref{corollary 3.18.1} we have  $\dashnorm w\|_{L_{p,q}(C_{ \rho_{2}} )}\leq N\rho_{2}^{ 2}\dashnorm \cL_{0} u\|_{L_{p,q}(C_{\rho_{2}})}$.

It follows  by H\"older's inequality that 
$$
\dashnorm D^{2}v\|_{L_{p,q}(C_{2\rho_{1}})}
\leq \dashnorm D^{2}v\|_{L_{r,r}(C_{2\rho_{1}})}
\leq  N   ( \rho_{2} /  \rho_{1})  ^{ (d+2)/r}\dashnorm D^{2}v\|_{L_{r,r}(C_{ \rho_{2}/2})}
$$
\begin{equation}
                              \label{7.16.4}
\leq N  (\rho_{2}/\rho_{1})^{\alpha}
\big[\rho_{2} ^{-2}\dashnorm u-l\|_{L_{p,q}(C_{ \rho_{2}} )}+ \dashnorm \cL_{0} u\|_{L_{p,q}(C_{ \rho_{2}} )}\big] ,
\end{equation}
where the last inequality is due to $(d+2)/r\leq\alpha$.

Then according to  
\eqref{10.4.4}  that
\begin{equation}
                              \label{7.18.1}
\dashnorm D^{2}u\|_{L_{p,q}(C_{\rho_{1}})}\leq
N\dashnorm \cL_{0} u\|_{L_{p,q}(C_{2\rho_{1}})}+
N \rho_{1}^{-2}\dashnorm u-l'\|_{L_{p,q}(C_{2\rho_{1}})},
\end{equation}
where $l'$ is any affine function.
Here $u-l'=w+(v-l')$, and the last term in \eqref{7.18.1} is dominated by   
$$
N(\rho_{2}/\rho_{1})^{ 2+d/p+2/q}\dashnorm \cL_{0} u\|_{L_{p,q}(C_{\rho_{2}})}+
N \rho_{1}^{-2}\dashnorm v-l'\|_{L_{p,q}(C_{2\rho_{1}})},
$$
where the second norm, for an appropriate choice of $l'$, is estimated by  a constant times 
$$
\dashnorm D^{2}v\|_{L_{p,q}(C_{2\rho_{1}})}+
\dashnorm \partial_{t}v\|_{L_{p,q}(C_{2\rho_{1}})}\leq N\dashnorm D^{2}v\|_{L_{p,q}(C_{2\rho_{1}})}
$$
owing to the Poincar\'e inequality.
After that, to get \eqref{7.15.3},
it only remains to refer to \eqref{7.16.4}.
Estimate \eqref{7.15.30} is obtained
from \eqref{7.15.3} by the Poincar\'e inequality.
The lemma is proved. \qed

\begin{lemma}
                     \label{lemma 7.16.1}
Let $p,q\in(1,\infty)$, $\beta>\alpha$,
  $\rho \leq 1 $,
$u\in W^{1,2}_{p,q }(C_{\rho  }) $ and set
$f:=\cL_{0}u $.
 Then
there exists $\nu=\nu(\alpha,\beta,d,\delta,p,q, \rho_{a})>1$
such that
for any  $r\leq\rho /\nu$
we have
\begin{equation}
                                 \label{7.18.2}
r^{\beta}\dashnorm \partial_{t}u, D^{2}u\|_{ L_{p,q}(C_{r})}
\leq N \sup_{r\leq s\leq\rho }
s^{\beta}\dashnorm f\|_{ L_{p,q}(C_{s })}
+ N\rho ^{\beta-2}\dashnorm u-l\|_{L_{p,q}(C_{\rho })},
\end{equation}
where $l$ is any affine function and the constants
 $N$ depend only on $\alpha$, $\beta$, 
$ \rho_{a}$,
$d,\delta,p,q$.
\end{lemma}

Proof. Since $\partial_{t}u=f-a^{ij}D_{ij}u$,
it suffices to concentrate on $D^{2}u$. Take the smallest  $\nu\geq 4$ such that
$N_{2}  \nu^{\alpha-\beta }\leq 1/2$ ($N_{2}$ is from \eqref{7.15.30}).
Then, for $r\leq \rho/\nu$, define 
$r_{n}=\nu^{n}r$, $m=\max\{n\geq 0:
r_{n+1}\leq\rho \}$,
$$
A_{n}=r_{n}^{\beta}\dashnorm D^{2}u\|_{ L_{p,q}(C_{r_{n}})},\quad B =\sup_{r\leq s\leq\rho }
s ^{\beta}\dashnorm f\|_{ L_{p,q}(C_{s})}.
$$
For $0\leq n\leq m$, $\rho_{1}=r_{n}$,
$\rho_{2}=r_{n+1}$ estimate \eqref{7.15.30}
yields
$$
A_{n }\leq N_{1}\nu^{\gamma-\beta}B +N_{2}\nu^{\alpha-\beta}
A_{n+1}\leq N_{3}B +(1/2)
A_{n+1}.
$$
By iterating we obtain
$A_{0}\leq 2N_{3}B+2^{- m}A_{m }$ and arrive at
\begin{equation}
                              \label{7.22.1}
r^{\beta}\dashnorm D^{2}u\|_{ L_{p,q}(C_{r})}
\leq N \sup_{r\leq s\leq\rho }
s ^{\beta}\dashnorm f\|_{ L_{p,q}(C_{s })}
+2^{- m}r_{m}^{\beta}\dashnorm
D^{2}u\|_{ L_{p,q}(C_{r_{m}})}.
\end{equation}
Here by \eqref{7.15.3}  
$$
r_{m}^{\beta}\dashnorm
D^{2}u\|_{ L_{p,q}(C_{r_{m}})}\leq N
r_{m+1}^{\beta}\dashnorm f\|_{ L_{p,q}(C_{r_{m+1}})}
+
Nr_{m+1}^{\beta-2}\dashnorm  u-l\|_{L_{p,q}(C_{r_{m+1}})}
$$
$$
\leq N \sup_{r\leq s\leq\rho }
s ^{\beta}\dashnorm f\|_{ L_{p,q}(C_{s })}
+N\rho^{\beta-2}\dashnorm  u-l\|_{L_{p,q}(C_{\rho})},
$$
where the last inequality is due to the fact that
$r\leq r_{m+1}\leq\rho$, $r_{m+1}/\rho\geq\nu^{-1}$.  This proves the lemma. \qed

Lemma \ref{lemma 7.16.1} (with $\rho=1$) and Remark \ref{remark 10.3.1}
immediately lead to the following
(recall that  
Assumption \ref{assumption 3.25.1} is supposed to be satisfied).

\begin{corollary}
                            \label{corollary 10.27.1}
Let $p,q\in(1,\infty)$  
and let $\beta>\alpha $.   Let 
$u\in W^{1,2}_{p,q,\loc} $. 
 Then
\begin{equation}
                                 \label{7.18.21}
\|\partial_{t}u, D^{2}u\|_{E_{p,q,\beta} }
\leq N \| \cL_{0} u\|_{E_{p,q,\beta} }
+ N \sup_{C\in\bC_{1}}\|u\|_{L_{p,q}(C)},
\end{equation}
where  the constants
 $N$ depend only on $\alpha$, $\beta$, 
$\rho_{a}$,
$d,\delta,p,q$.
\end{corollary}

Next step is to replace $\cL_{0}u$ in \eqref{7.18.21}
and drop the last term.
For that we need the following where $d+1$ in place of $d$
comes into play.

\begin{lemma}
                           \label{lemma3.26.1}
Let   $d/p+2/q\geq\beta>\alpha$.
Then there are 
$$
(\tilde \lambda_{0},\tilde N_{0})=(\tilde \lambda_{0},\tilde N_{0})(\alpha,\beta,d,\delta,p,q, \rho_{a}),\quad \tilde \lambda_{0}>0,
$$
such that for any $\lambda\geq \tilde\lambda_{0}$ and any
$u\in W^{1,2}_{p,q,\loc} $, such that
\begin{equation}
                                 \label{3.26.5}
\sup_{C\in\bC_{1}}\|u\|_{L_{p,q}(C)}<\infty,
\end{equation}
we have 
\begin{equation}
                                 \label{6.3.02}
\|\lambda u,\sqrt\lambda Du, D^{2}u\|_{E_{p,q,\beta} }\leq \tilde N_{0}
\|\cL_{0} u -\lambda u\|_{E_{p,q,\beta}}.
\end{equation}

\end{lemma}

Proof. We use the same Agmon's method as in the proof of 
Lemma 6.3.8 of \cite{Kr_08}.
Consider the space
$$
\bR^{d+2}=\{(t,z)=(t,x,y):t,y\in\bR,x\in\bR^{d}\}
$$
 and the function
$$
\tilde{u}(t,z)=u(t,x) \cos(\mu y),
$$
 where
 $\mu=\sqrt{\lambda}$, $\lambda\geq0$.
Also introduce the operator
$$
\tilde{\cL}_{0}v(t,z)=\cL_{0}(t,x)v(t,z)+D_{yy}v (t,z).
 $$ 

Obviously the above theory is applicable to
$\tilde{u}$ and $\tilde{\cL_{0}}$ and we get
from \eqref{7.18.21}
$$
\|D_{yy}\tilde{u} \|_{E_{p,q,\beta}(\bR^{d+2}) }
\leq N \|\tilde{L_{0}}\tilde{u}  \|_{E_{p,q,\beta}(\bR^{d+2}) }
+N \sup_{C\in\bC_{1}}\|u\|_{L_{p,q}(C)}.
$$

Observe that for any $\rho<1, C\in\bC_{\rho}$
and $y\in \bR$
$$
\dashnorm u\|_{L_{p,q}(C)} \Big(\dashint_{y}^{\,\,y+\rho}|\cos(\mu r)|^{p}\,dr\Big)^{1/p}
\geq \gamma \dashnorm u\|_{L_{p,q}(C)},
$$
where $\gamma>0$ is an absolute constant. It follows that 
$$
\lambda \|u \|_{E_{p,q,\beta}  }
\leq N \|\tilde{L_{0}}\tilde{u}  \|_{E_{p,q,\beta}(\bR^{d+2}) }
+N \sup_{C\in\bC_{1}}\|u\|_{L_{p,q}(C)}.
$$
Since
$
\tilde{\cL_{0}}\tilde{u}  = \cos(\mu y)
[\cL_{0}u-\lambda u]$,
we have
$$
\|\tilde{L_{0}}\tilde{u}  \|_{E_{p,q,\beta}(\bR^{d+2})}\leq
N\|\cL_{0} u -\lambda u\|_{E_{p,q,\beta}},
$$
so that
\begin{equation}
                                                     \label{06.6.27.9}
\lambda \|u \|_{E_{p,q,\beta}  }
\leq N \|\cL_{0}u- \lambda u  \|_{E_{p,q,\beta} }
+N_{1} \sup_{C\in\bC_{1}}\|u\|_{L_{p,q}(C)}.
\end{equation} 
For $\lambda\geq 2N_{1}$ we can absorb the last term in \eqref{06.6.27.9} (which is finite)  into the 
left-hand side
and get the claimed estimate of $\lambda u$.
Furthermore, then \eqref{7.18.21} implies that
$$
\|D^{2}u\|_{L_{p,q,\beta}}\leq 
N \| \cL_{0} u-\lambda u\|_{E_{p,q,\beta} }
+N\lambda \|u\|_{E_{p,q,\beta}}  
+ N \sup_{C\in\bC_{1}}\|u\|_{L_{p,q}(C)},
$$
which yields the estimate of $D^{2}u$.
After that it only remains to use 
  Lemma \ref{lemma 3.27.4}.
The lemma is proved. \qed

Before the following existence theorem it is useful
to note that the functions $u$ in $ E^{1,2}_{p,q,\beta}$
are also in 
$  W^{1,2}_{p,q,\loc}$ and \eqref{3.26.5}
holds for them.
\begin{theorem}
                         \label{theorem 3.27.1}
Let   $d/p+2/q\geq\beta>\alpha$, $\lambda\geq
\tilde \lambda_{0}$. Then for any $f\in E_{p,q,\beta}$
there is a unique $u\in E^{1,2}_{p,q,\beta}$
  such
that $\cL_{0}u-\lambda u=f$.

\end{theorem}

Proof. To prove uniqueness, assume that $u\in E^{1,2}_{p,q}$
and $\cL_{0}u-  \lambda u=0$. Then take $\zeta\in
C^{\infty}_{0}(\bR^{d+1})$ such that $0\leq \zeta\leq \zeta(0)=1$
and set $\zeta_{n}(x)=\zeta(x/n)$, $u_{n}=u\zeta_{n}$.
We have that $u_{n}\in W^{1,2}_{p,q}$ and
$$
|\cL_{0}u_{n}-\lambda u_{n}|\leq N(n^{-1}|Du|+
n^{-2}|u|),
$$
where $N$ is independent of $n$.
By Lemma  \ref{lemma3.26.1}, we have $\|u_{n}\|_{L_{p,q,\beta}}
\leq N(n^{-1}+n^{-2})$ and letting $n\to
\infty$ proves the uniqueness of $E^{1,2}_{p,q}$-solutions.

To prove the existence part,
take $f\in E_{p,q,\beta}$ and set $f_{n}=f\zeta_{n}$.
Then by Theorem \ref{theorem 12.18.1}, for each $n$,
there exists (a unique) $u_{n}\in W^{1,2}_{p,q}$
such that $\cL_{0}u_{n}-\lambda u_{n}=f_{n}$.
By Lemma  \ref{lemma3.26.1}  we have
that $u_{n}\in E^{1,2}_{p,q,\beta}$ and  
 
\begin{equation}
                                \label{3.20.02}
\|\partial_{t} u_{n},D^{2}u_{n},  Du_{n},  u_{n}\|_{E_{p,q,\beta }}
\leq N\|f\zeta_{n}\|_{E_{p,q,\beta }}
\leq N\|f \|_{E_{p,q,\beta }},  
\end{equation}
which implies that there is a subsequence $n'\to\infty$ such that in any $C\in\bC$
weakly in $L_{p,q}(C)$
\begin{equation}
                                \label{3.20.03}
(\partial_{t} u_{n'},D^{2}u_{n'},  Du_{n'},  u_{n'})
\to (\partial_{t} u ,D^{2}u ,  Du ,  u ),
\end{equation}
where $u\in W^{1,2}_{p,q,\loc}$. Then \eqref{3.20.02}
also implies that $u\in E^{1,2}_{p,q,\beta}$
and along with \eqref{3.20.03} and the fact that  $f_{n}\to f$
in any $L_{p,q}(C)$ for any   $C\in \bC$
show that $\cL_{0}u-\lambda u=f$. 
The theorem is proved. \qed

 \mysection{Proof of Theorem \protect\ref{theorem 10.27.1}}
                         \label{section 4.2.3}

  Here we return to the set up in \eqref{3.21.01}
and assumptions introduced in Section \ref{section 3.27.3}.
Recall that $\rho_{b,c}=1$.
To deal with the full operator $\cL$ we need the following, in which   we need to distinguish the cases $\beta_{0}>1$ and $\beta_{0}<1$ for the first time in this article dealing with the  equations and recall that
$\beta_{0}\leq d/p+2/q$.

\begin{lemma}
                          \label{lemma 10.26.1}  
For any $u\in E^{1,2}_{p,q,\beta_{0}}$ we have
for $\beta_{0}>1$ that
\begin{equation}
                         \label{2.29.1}
\| b|Du|\,\|_{E_{p,q,\beta_{0}}}
 \leq N\mathfrak{b}_{1}
\|\partial_{t}u,D^{2} u\|_{E_{p,q,\beta _{0}}}
+N\mathfrak{b}_{1}  
\| u\|_{E_{p,q,\beta_{0} }},
\end{equation}
where $N=N(d,p ,q ,\beta_{0} )$,
and for $\beta_{0}<1$ and  any $\varepsilon
\in(0,1]$
\begin{equation}
                                \label{4.1.7}
\| b|Du|\,\|_{E_{p,q,\beta_{0}}}
\leq \varepsilon^{1-\beta}
\|\partial_{t}u,D^{2} u\|_{E_{p,q,\beta _{0}}}
+N\varepsilon^{-1}\mathfrak{b}^{1+\beta/(1-\beta)}\| u\|_{E_{p,q,\beta_{0} }},
\end{equation}
where  
$N=N(d,p,q,\beta_{0})$.
\end{lemma}

Proof. To prove the first assertion let $\zeta=I_{C_{1}}$ and $\xi\in C^{\infty}_{0}(\bR^{d+1})$ be such that $\xi=1$ on $C_{1}$.
It is easy to see
that for any $\rho>0$ and $C\in\bC_{\rho}$
(here we use that $\beta_{0}\leq d/p +2/q $ when $\rho\geq 1$) 
\begin{equation}
                      \label{10.27.2}
\dashnorm \zeta b\|_{L_{p\beta_{0},q\beta_{0}}(C)}
\leq N(d )\mathfrak{b}_{1} \rho^{-1}.
\end{equation}
Therefore, by Remark 5.8 of \cite{Kr_22},
since $1<\beta_{0}\leq d/p+2/q$, 
\begin{equation}
                      \label{10.27.1}
\| \zeta b|D(\xi u)|\,\|_{E_{p,q,\beta_{0}}}
 \leq N\mathfrak{b}_{1}\sup_{\rho>0}\sup_{C\in\bC_{\rho}}
\rho^{\beta_{0}}\dashnorm(\partial_{t}+\Delta)(\xi u)\|_{L_{p,q}(C)}.
\end{equation}
Similarly to \eqref{10.27.2} the sup in the right-hand side
is dominated by
$$
N\|(\partial_{t}+\Delta)(\xi u)\|_{E_{p,q,\beta_{0}}}
\leq N\|(\partial_{t}+\Delta)u\|_{E_{p,q,\beta_{0}}}
+N \| u\|_{E_{p,q,\beta_{0}}}
+N \| Du\|_{E_{p,q,\beta_{0}}}.
$$

By combining this with \eqref{3.20.6},
for   $\varepsilon=1$, we obtain
\eqref{2.29.1} with the left-hand side
replaced by the left-hand side of 
\eqref{10.27.1}. Thus obtain inequality,
obviously holds if we take any translation of $\zeta,\xi$.
After that it only remains to observe that
the supremum of the left-hand sides of 
\eqref{10.27.1} over translated $\zeta,\xi$
is definitely larger than the left-hand side of
\eqref{2.29.1}. 

The second assertion follows directly from Lemma
\ref{lemma 3.27.40} (B).
The lemma is proved. \qed

\begin{lemma}
                     \label{lemma 3.29.1}
For any $u\in E^{1,2}_{p,q,\beta_{0}}$
we have for $\beta_{0}>2$ that
$$
\|cu\|_{E_{p,q,\beta_{0}}}\leq N(d,p,q,\beta_{0})
\mathfrak{c}_{1}\|u\|_{E^{1,2}_{p,q,\eta_{0}}},
$$
and for $\beta_{0}< 2$ that for any
  $\varepsilon\in(0,1]$
\begin{equation}
                                  \label{3.29.1}
\|cu\|_{E_{p,q,\beta_{0}}}\leq \varepsilon
\|\partial_{t}u,D^{2}u\|_{E_{p,q,\beta}}
+N \nu\| u\|_{E_{p,q,\beta}}
\end{equation}
where
$N=N( d, p,q,\beta_{0})$, $\nu= \nu(\mathfrak{c} ,\varepsilon,\beta_{0})=\mathfrak{c}^{2/(2-\beta_{0})}\varepsilon^{-\beta_{0}/(2-\beta_{0})}$.
 
\end{lemma}

Proof. First assume $\beta_{0}>2$.
Set $(r,s)=\beta_{0}(\beta_{0}-2)^{-1}(p  ,q  )$. Then by H\"older's inequality
$$
\rho^{\beta_{0}}\dashnorm c  u\|_{L_{p,q }(C_{\rho})}
\leq \rho^{2}\dashnorm c  \|_{L_{p\beta_{0},q\beta_{0}}(C_{\rho})}\cdot \rho^{\beta_{0}-2}
\dashnorm  u\|_{L_{r ,s}(C_{\rho})}, 
$$
so that
$$
\|cu\|_{E_{p,q,\beta_{0}}}\leq \mathfrak{c}
\|u\|_{E_{r,s,\beta_{0}-2}},
$$
and it only remains to apply  Lemma \ref{lemma 4.1.1}.

  In case $\beta_{0}<2$
our assertion follows directly from
Lemma \ref{lemma 3.27.40} (A). The lemma is proved. \qed

{\bf Proof of Theorem \ref{theorem 10.27.1}}. In light of Theorem \ref{theorem 3.27.1} 
and the method of continuity, it suffices to prove 
\eqref{10.28.01} as an a priori estimate.
We will show how the proof works only when $\beta_{0}>2$.
In the other cases: $1<\beta_{0}<2$ and $\beta_{0}<1$,
the proof is quite similar.

By using Lemmas \ref{lemma3.26.1},   \ref{lemma 10.26.1}, and \ref{lemma 3.29.1},
we see  that for $\lambda\geq \tilde \lambda_{0}$
and $\varepsilon\in(0,1]$
$$
\|\partial_{t} u,D^{2}u, \sqrt\lambda Du, \lambda u\|_{E_{p,q,\beta_{0} }}
\leq \tilde N_{0}\|\cL u- \lambda  u\|_{E_{p,q,\beta_{0} }}+\tilde N_{0}\|b^{i}D_{i}u\|_{E_{p,q,\beta_{0} }}
$$
$$
+\tilde N_{0}\|c u\|_{E_{p,q,\beta_{0} }}
\leq \tilde N_{0}\|f\|_{E_{p,q,\beta_{0} }}+
\tilde N_{0} N_{1}(\mathfrak{b}_{1}
+\mathfrak{c}_{1})\big( \|\partial_{t}u,D^{2}u \|_{E_{p,q.\beta_{0}}}
+\|u\|_{E_{p,q.\beta_{0}}}\big)
$$
where $ N_{1} = N_{1} (d,p ,q , \beta_{0} )$.
It follows that, if $\tilde N_{0}N_{1}(\check b+\check c)
\leq 1/2$, then for
 $\mathfrak{b}_{1}\leq \check b$, $\mathfrak{c}_{1}
\leq \check c$, and
$
\lambda\geq 1
$
we get \eqref{10.28.01}. The theorem is proved. \qed

\mysection{Proof of Lemmas  \protect\ref{lemma 3.27.40} and \protect\ref{lemmq 4.6.1}}

                           \label{section 4.2.1}

For $s,r>0,\alpha>0$, and appropriate $f(t,x)$'s
on $\bR^{d+1}$ define
$$
p_{\alpha}(s,r)=\frac{1}{s^{(d+2-\alpha)/2}}e^{-r^{2}/s}I_{s>0},
$$
$$
P_{\alpha}f(t,x)=\int_{\bR^{d+1} }p_{\alpha}(s,|y|)f(t+s,x+y)\,dyds.
$$

{\bf Proof of Lemma \ref{lemma 3.27.40}}.
Case  (A), (a). Take $\zeta\in C^{\infty}_{0}$
such that $\zeta=1$ in $C_{1/2}$, $\zeta(t,x)=0$
in $\{t\geq0\}\cap C_{1}^{c}$, and $\zeta$ depends only
on $t$ and $|x|$.  For $\varepsilon
\in(0,1]$ define $\zeta_{\varepsilon}(t,x)=\zeta(t/\varepsilon^{2},x/\varepsilon)$. Introduce
  $f=-\partial_{t}u-(1/4)\Delta u$. Then
   $u(0)=u^{\varepsilon}_{1} 
+u_{2}^{\varepsilon} $, where with a constant $\gamma=
\gamma(d)$ and $\bR^{d+1}_{+}=(0,\infty)\times\bR^{d}$
$$
\gamma^{-1}u^{\varepsilon}_{1}=
\int_{\bR^{d+1}_{+} }\zeta_{\varepsilon}(t,x)t^{-d/2}
e^{-|x|^{2}/t}f(t ,x )\,dxdt,
$$
$$
\gamma^{-1}u_{2}^{\varepsilon} =
\int_{\bR^{d+1}_{+} }(1-\zeta_{\varepsilon}(t,x))t^{-d/2}
e^{-|x|^{2}/t}f(t ,x )\,dyds
$$
$$
=\int_{\bR^{d+1}_{+} }u(t ,x )\big(\partial_{t}-(1/4)\Delta )\big[(1-\zeta_{\varepsilon}(t,x))t^{-d/2}
e^{-|x|^{2}/t}\big]\,dxdt
$$
\begin{equation}
                                  \label{3.31.2}
=  \int_{\bR^{d+1}_{+} }u(t ,x ) \eta _{\varepsilon}(t,x) t^{-d/2}
e^{-|x|^{2}/t}\,dxdt,
\end{equation}
where
$$
\eta_{\varepsilon}(t,x)= -(\partial_{t}-(1/4)\Delta ) \zeta_{\varepsilon}(t,x))-x^{i}D_{i}\zeta_{\varepsilon}(t,x) t^{-1}.
$$
Note that $\eta_{\varepsilon}=0$ outside the part of $C_{\varepsilon}$ lying in $\bR^{d+1}_{+}$, $D\zeta(t,0)=0$, 
and 
$$
|x^{i}D_{i}\zeta_{\varepsilon}(t,x)| t^{-1}
\leq N\varepsilon^{-2}|x|^{2}/t,\quad e^{-|x|^{2}/t}|x|^{2}/t
\leq Ne^{-|x|^{2}/(2t)}.
$$

Therefore,
by inspecting
the proof of Lemma 2.2 of \cite{Kr_22}
dealing with  $P_{\alpha}(fI_{C_{\rho}})$
when $\alpha=2$, $\rho=\varepsilon$,   one easily sees
that
$$
\gamma^{-1}|u(0)|\leq P_{2}( |f\zeta_{\varepsilon}|I_{C_{\varepsilon}})
+N P_{2}( |u\eta _{\varepsilon}|I_{C_{\varepsilon}})\leq N\varepsilon^{2}\dashint_{C_{\varepsilon}}
|f|\,dxdt
$$
$$
+N\int_{0}^{\varepsilon}r\dashint_{C_{r}}|f|\,dxdt\,dr
+N \dashint_{C_{\varepsilon}}
|u|\,dxdt+N\varepsilon^{-2}\int_{0}^{\varepsilon}r\dashint_{C_{r}}|u|\,dxdt\,dr
$$
$$
\leq N\|f\|_{E_{p,q,\beta}}\big(\varepsilon^{2-\beta}+
\int_{0}^{\varepsilon}r^{1-\beta}\,dr\Big)
+N\varepsilon^{-\beta}\|u\|_{E_{p,q,\beta}}.
$$
This yields \eqref{3.20.06} at the origin
and, certainly, proves it everywhere. To prove the continuity of $u$ it suffices to observe that
$u^{(\varepsilon)}\to u$ uniformly in light
of \eqref{3.20.6} and Lemma \ref{lemma 3.14.3},
in which we can take any $\beta'>\beta,\beta'<2$.
This proves
the lemma in Case (A), (a).

In Case (A), (b) we have $u\in W^{1,2}_{p,q}(C_{1})$
and by the well-known embedding theorems we have \eqref{3.20.06} with $L_{p,q}$ in place of $E_{p,q,\beta}$, which yields \eqref{3.20.06} as is since
$$
 \|u\|_{L _{p,q}(C_{1})}\leq N(d)
\|u\|_{E {p,q,\beta}}.
$$

Case (B), (a). Here we repeat the first part
of the proof with $De^{-|x|^{2}/t}=-(2x/t)
e^{-|x|^{2}/t}$ in place of $e^{-|x|^{2}/t}$.
This   leads to an analog of \eqref{3.31.2} in which
we have to replace $\eta_{\varepsilon}$
with
$$
\bar\eta_{\varepsilon}(t,x)=-(2x/t)\eta_{\varepsilon} 
-2(1-\zeta_{\varepsilon})x/t^{2} -4D_{i}\zeta_{\varepsilon}D_{i}x/t.
$$
Again $\xi_{\varepsilon}=0$ outside the part of $C_{\varepsilon}$ lying in $\bR^{d+1}_{+}$ and
$$
|\xi_{\varepsilon}|\leq N\varepsilon^{-2}(|x|/t^{1/2}+
|x|^{3}/t^{-3/2})t^{-1/2},
$$
$$
(|x|/t^{1/2}+
|x|^{3}/t^{-3/2})t^{-1/2}t^{-d/2}e^{-|x|^{2}/t}
\leq Nt^{-(d+1)/2}e^{-|x|^{2}/(2t)}.
$$
This time using
the proof of Lemma 2.2 of \cite{Kr_22}
with $\alpha=1$ and denoting $\bar\zeta_{\varepsilon}
(t,x)=-2\zeta_{\varepsilon}
(t,x)x/t^{1/2}$ we obtain
$$
\gamma^{-1}|Du(0)|\leq P_{1}( |f\bar\zeta_{\varepsilon}|I_{C_{\varepsilon}})
+N\varepsilon^{-2}P_{1}( |u\bar\eta_{\varepsilon}|I_{C_{\varepsilon}})\leq N\varepsilon\dashint_{C_{\varepsilon}}
|f|\,dxdt
$$
$$
+N\int_{0}^{\varepsilon} \dashint_{C_{r}}|f|\,dxdt\,dr
+N \varepsilon^{-1}\dashint_{C_{\varepsilon}}
|u|\,dxdt+N\varepsilon^{-2}\int_{0}^{\varepsilon} \dashint_{C_{r}}|u|\,dxdt\,dr
$$
$$
\leq N\|f\|_{E_{p,q,\beta}}\big(\varepsilon^{1-\beta}+
\int_{0}^{\varepsilon}r^{-\beta}\,dr\Big)
+N\varepsilon^{-\beta-1}\|u\|_{E_{p,q,\beta}}.
$$
This yields \eqref{3.20.60} in Case (B), (a).
The continuity of $Du$ is derived as in Case (A), (a).

In Case (B), (b) our assertion again follows from embedding theorems for $W^{1,2}_{p,q}$.
The lemma is proved.  \qed 

{\bf Proof of Lemma \ref{lemmq 4.6.1}}. Use the notation $(f)_{r}$ for the average of $f$ over $C_{r}$
and observe that by Poincar\'e's inequality
$$
\dashint_{C_{r}}|u-(u)_{r}-x^{i}(D_{i}u)_{r}|
\,dxdt\leq N r^{2}\dashint_{C_{r}}|\partial_{t}u|\,dxdt+Nr\dashint_{C_{r}}
|Du-(Du)_{r}|\,dxdt,
$$
\begin{equation}
                        \label{4.6.3}
\dashint_{C_{r}}
|Du-(Du)_{r}|\,dxdt\leq Nr\dashint_{C_{r}}(
|\partial_{t}u|+|D^{2}u|)\,dxdt.
\end{equation}
Also by Lemma \ref{lemma 3.27.5} and H\"older's
inequality for $r\leq 1$
$$
\dashint_{C_{r}}| x^{i}(D_{i}u)_{r}|
\,dxdt\leq Nr\dashint_{C_{r}}|Du|\,dxdt
\leq Nr^{2-\beta}\|u\|_{E^{1,2}_{p,q,\beta}}.
$$
By combining these estimates we see that for $r\leq1$
$$
\dashint_{C_{r}}|u-(u)_{r} |
\,dxdt\leq Nr^{2-\beta}\|u\|_{E^{1,2}_{p,q,\beta}}.
$$
This estimate also holds for translates of $C_{r}$.
After that, to prove \eqref{4.6.1}, it only remains to invoke Campanato's
result. His result also proves \eqref{4.6.2}
owing to \eqref{4.6.3}.
The lemma is proved. \qed

{\bf Acknowledgements}. The author is sincerely grateful to Hongjie Dong who pointed out the possibility to
use Morrey spaces when the drift term has poor integrability.
\medskip

The paper has no  data and there is no conflict of interest.

\end{document}